# PERIODIC ORBITS, SYMBOLIC DYNAMICS AND TOPOLOGICAL ENTROPY FOR THE RESTRICTED 3-BODY PROBLEM

GIANNI ARIOLI

Abstract. This paper concerns the restricted 3-body problem. By applying topological methods we give a computer assisted proof of the existence of some classes of periodic orbits, the existence of symbolic dynamics and we give a rigorous lower estimate for the topological entropy.

## 1. Introduction

The problem of both restricted and full $N$-body systems has such a long story that it is impossible to give an extensive bibliography here; we refer the reader to the classical texts [MH, M, SM, S].

This paper concerns the study of periodic and chaotic solutions of the planar restricted 3-body problem. If we assume that the primaries orbit around each other with period $2\pi$ and we use a rotating reference frame, i.e. if we use synodical coordinates, then the motion of the third body is described by the following system of second order differential equations:

$$\begin{aligned} \ddot{x} + 2\dot{y} &= \Omega_x \\ \ddot{y} - 2\dot{x} &= \Omega_y, \end{aligned} \tag{1.1}$$

where

$$\Omega(x,y) = \frac{x^2}{2} + \frac{y^2}{2} + \frac{m_1}{\sqrt{(x-R_1)^2+y^2}} + \frac{m_2}{\sqrt{(x+R_2)^2+y^2}} + C,$$

$m_1$ and $m_2$ are the masses of the primaries, $C$ is an arbitrary constant and $R_1$, $R_2$ depend on $m_1$ and $m_2$.

The methods we employ are topological with computer assistance. This kind of computer assisted proof has been introduced in [MM1, MM2, MM3], while the topological methods employed here have been introduced in [Z1, Z2], see also the applications in [AZ, GZ, Z3]. We point out that, although some of the proofs we give require computer assistance, they are rigorous and they can be easily reproduced on any recent computer. To this purpose, in the last section we provide all necessary information, mostly by referring to other papers, and we also provide on the web a Mathematica version of the algorithms employed.

The main result presented in this paper is the proof of symbolic dynamics on fourteen symbols which corresponds to the existence of orbits, periodic or non-periodic, that come close to five different periodic orbits in any prescribed order. This result yields a lower estimate for the topological entropy of the system, and it is therefore a proof of its chaotic behavior. To the author knowledge, no rigorous estimate of the topological entropy for the Poincaré map of this system is available in the literature. This paper serves other purposes as well. We show how to extend the computer assisted techniques developed in the papers cited above to the planar restricted 3-body problem, which presents different features and difficulties and requires

This research was supported by MURST project "Metodi variazionali ed Equazioni Differenziali Non Lineari".

new techniques; to this purpose we introduce some new topological tools and computational methods. We give a rigorous proof of the existence of a class of periodic orbits at different energy levels: the existence of such periodic orbits (and many more) is well known, but to the author knowledge the proof is sometimes purely numerical, with no mathematical rigour, or perturbative, with no estimate on the range of validity of the perturbation parameter. We also provide a very narrow estimate of the locations of the intersection of such periodic orbits with the line connecting the primaries. In summary, this paper provides both the methods and some examples on how the computer assisted techniques can yield results for the planar restricted 3-body problem. We point out that these results are not of pertubative nature, and to the author knowledge the result on symbolic dynamic is not accessible by purely analytical methods; in particular the problem is not treatable by Melnikov's method.

A result on chaos for the planar restricted 3-body problem has been recently presented in [SK] with a different method, also based on computer assisted techniques, but the authors do not claim to give a rigorous proof. Indeed, they use, in their words, "realistic estimated upper bounds" for the errors made in the numerical computation of the Poincaré map. Such method found a rigorous application on a different system in [KMS]. Here we estimate rigorously all computational errors and we do provide rigorous proofs of all the theorems we state.

The paper is organized as follows: in Section 2 we provide a brief introduction to the problem; in Section 3 we introduce the Poincaré maps that are used in the proofs and explain a symmetry of the system which is widely employed in the proofs; in Section 4 we explain our method for detecting and proving periodic orbits and present the results we obtain. Section 5 is the main portion of the paper. Here we introduce the topological and computational methods and we describe the results on symbolic dynamics and on topological entropy. Ideas on further developments are given in Section 6 and details on the computer assisted proofs are in Section 7.

## 2. Description of the system

We first derive briefly equation (1.1). It is well known that the two body problem admits a solution where both bodies move in a circular counterclockwise motion around their center of mass (and of course it also admits the symmetric clockwise solution). We call the two bodies $P_1$ and $P_2$ (the primaries). If $P_1$ and $P_2$ have mass $m_1$ and $m_2$, then the circular solution with minimal period $2\pi$ minimizes the Lagrangian functional

$$f(x_1, x_2) := \int_0^{2\pi} \frac{m_1}{2}|\dot{x}_1|^2 + \frac{m_2}{2}|\dot{x}_2|^2 + \frac{m_1 m_2}{|x_1 - x_2|} dt$$

defined on $H^1_{per}([0, 2\pi])$, the space of $2\pi$-periodic functions whose weak derivative is square-integrable.

Take

$$x_1 = R_1(\cos t, -\sin t) \text{ and } x_2 = R_2(-\cos t, \sin t);$$

in order to compute $R_1$ and $R_2$ we have to minimize

$$\tilde{f}(R_1, R_2) = \frac{m_1}{2}R_1^2 + \frac{m_2}{2}R_2^2 + \frac{m_1 m_2}{R_1 + R_2}. \tag{2.1}$$

Setting 0 as the center of mass we have $R_2 = \frac{m_1}{m_2}R_1$, and by minimizing $\tilde{f}(R_1, R_2)$ we get $R_1 = \frac{m_2}{(m_1+m_2)^{2/3}}$ and $R_2 = \frac{m_1}{(m_1+m_2)^{2/3}}$.

The Lagrangian for the restricted 3-body problem is $\tilde{L}(w) = \frac{1}{2}|\dot{w}|^2 - V(w)$ where

$$V(w) = -\frac{m_1}{|w - x_1|} - \frac{m_2}{|w - x_2|}.$$

It is convenient to use synodical coordinates, i.e. a reference frame where the primaries sit still. Set $w = Rv$, where $R = \begin{pmatrix} \cos t & \sin t \\ -\sin t & \cos t \end{pmatrix}$ and let $J = \begin{pmatrix} 0 & -1 \\ 1 & 0 \end{pmatrix}$ be the standard symplectic matrix; then $\dot{w} = \dot{R}v + R\dot{v} = R(\dot{v} - Jv)$, therefore the Lagrangian is given by

$$L(v) = \frac{1}{2}|\dot{v}|^2 - (\dot{v}, Jv) + \Omega(v),$$

where

$$\Omega(v) = \frac{|v|^2}{2} + \frac{m_1}{|v - (R_1, 0)|} + \frac{m_2}{|v - (-R_2, 0)|} - C \tag{2.2}$$

It is well known that $H(x, \dot{x}, y, \dot{y}) = \dot{x}^2 + \dot{y}^2 - 2\Omega(x, y)$ is an integral of the motion (the Jacobi integral), therefore the motion takes place on the manifold $H(x, \dot{x}, y, \dot{y}) = h$. With some abuse of language, we call the Jacobi integral $H(v) = |\dot{v}|^2 - 2\Omega(v)$ the *energy*. The Euler-Lagrange equations are:

$$\ddot{v} - 2J\dot{v} - \nabla\Omega(v) = 0. \tag{2.3}$$

This system admits five equilibrium points. Three such points are aligned with the primaries and are called $L_1$, $L_2$ and $L_3$ (the collinear equilibrium solutions). One of the collinear points lies between the primaries; let that point be $L_1$. All the collinear points are saddles for $\Omega$. The remaining two equilibrium points are called $L_4$ and $L_5$; they are the absolute minima of $\Omega$ and they both sit at the vertices of an equilateral triangle with the primaries at the other vertices (the equilateral equilibrium solutions).

We only consider the case $m_1 = m_2$. The reason for this choice is that this case has been already widely considered and many families of periodic orbits for this case are known to exist. Nonetheless, rigorous estimates on the location of the periodic orbits, results on symbolic dynamics and topological entropy and even a rigorous proof of existence of the periodic orbits are new to the author knowledge. Furthermore the equal masses case is, in some sense, the opposite of a "perturbative" case with one primary much less massive than the other one. In this case $R_1 = R_2 = 2^{-2/3}$ and we choose the constant $C = -2^{5/3}$ in order to have $\Omega(L_1) = 0$. By a direct computation it is easy to see that there exists $h_0 = -2\Omega(L_2) = -2\Omega(L_3) \simeq .8623$ such that the region $\{v : 2\Omega(v) + h \geq 0\}$ where the motion can take place is split in exactly two connected regions, one bounded and the other unbounded, if and only if $0 \leq h < h_0$. We only consider values of $h$ in the range $[0, h_0)$ and we look for trajectories in the bounded region. Note that, if $h = 0$, the bounded region is equal to the union of two closed sets whose only intersection is the origin and no trajectory can intersect both sets; if $0 < h < h_0$ the bounded region is homeomorphic to a ball (if we include the singularities in the region); if $h < 0$ there are three connected regions and finally if $h \geq h_0$ the region $\{v : 2\Omega(v) + h \geq 0\}$ is connected, but not simply connected. The following pictures represent the curves $\{2\Omega(v) + h = 0\}$ for different values of $h$ (the small disks represent the primaries):

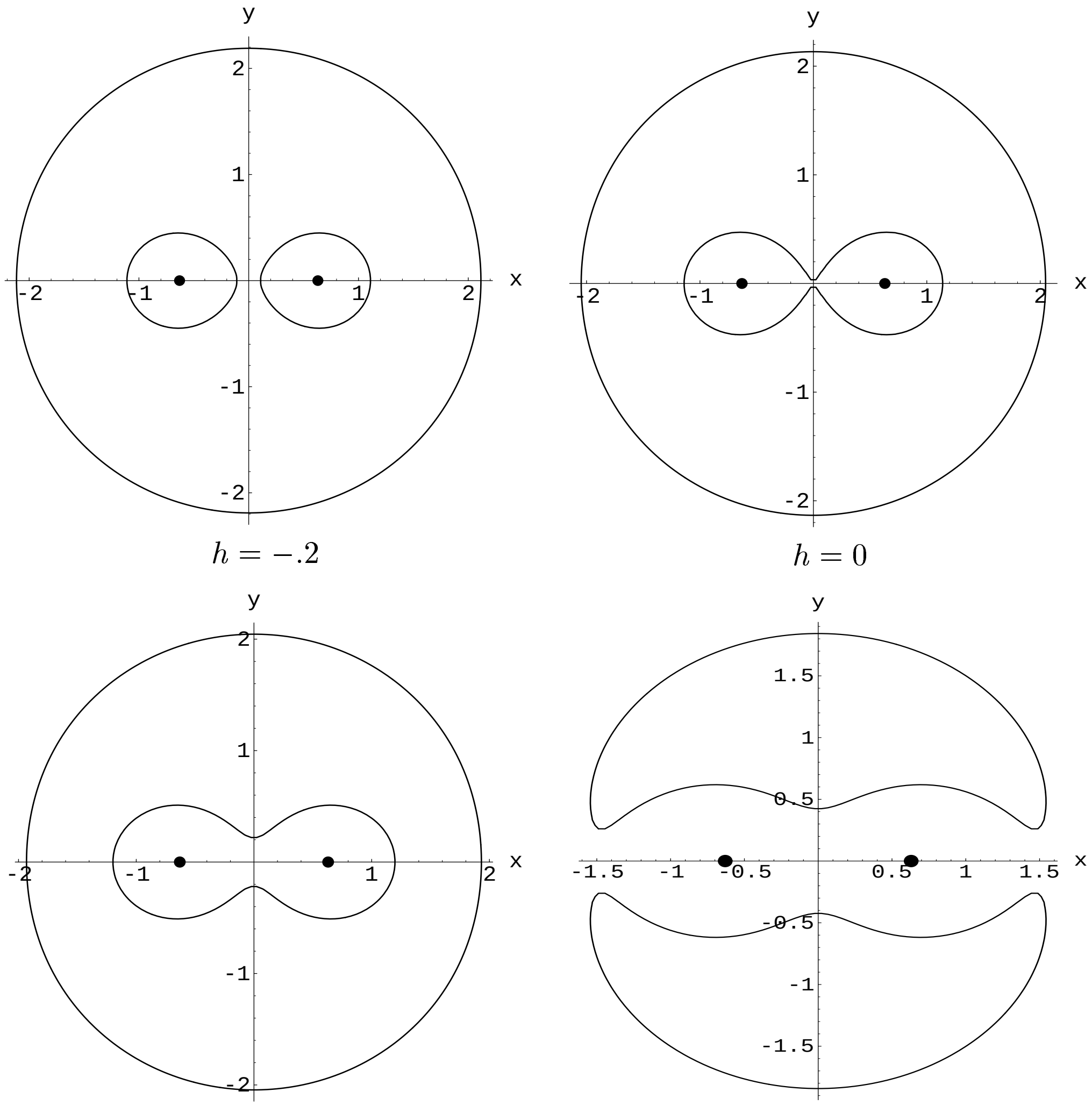

$h = -.2$ $h = 0$

$h = .3$ $h = .9$

## 3. The Poincaré maps

In order to study the system at some fixed energy $h$ we consider the Poincaré return map $P : D(P) \subset \mathbb{R}^2 \to \mathbb{R}^2$ defined in the following way. Given $(x, p_x)$ such that $x \neq \{R_{1,} - R_2\}$ and $2\Omega(x,0) + h - p_x^2 > 0$ there exists a unique positive value of $p_y = p_y(x, p_x) = \sqrt{2\Omega(x,0) + h - p_x^2}$ such that $H(x, p_x, 0, p_y) = h$. Let $\varphi(x, p_x; \cdot) : \mathbb{R} \to \mathbb{R}^4$ be the solution of the equation (1.1) with initial conditions $x(0) = x$, $\dot{x}(0) = p_x$, $y(0) = 0$, $\dot{y}(0) = p_y$ and call $\varphi_i$, $i = 1, \ldots, 4$ its components. By definition $\varphi_3(x, p_x; 0) = 0$ and $\varphi_4(x, p_x; 0) > 0$, therefore $\varphi_3(x, p_x; t) > 0$ for all positive and small $t$. If there exists a time $T_1$ such that $\varphi_3(x, p_x; T_1) = 0$ and $\varphi_3(x, p_x; t) > 0$ for all $t \in (0, T_1)$, and a time $T_2 > T_1$ such that $\varphi_3(x, p_x; T_2) = 0$ and $\varphi_3(x, p_x; t) < 0$ for all $t \in (T_1, T_2)$, then we define

$$P(x, p_x) = (\varphi_1(x, p_x; T_2), \varphi_2(x, p_x; T_2)).$$

In other words, $P$ is the return map on the section $y = 0$, $p_y > 0$.

In fact, for different reasons which we point out later, we find it useful to consider the half Poincaré maps $H_1 : D(H_1) \subset \mathbb{R}^2 \to \mathbb{R}^2$ and $H_2 : D(H_2) \subset \mathbb{R}^2 \to \mathbb{R}^2$. $H_1$ is defined by

$$H_1(x, p_x) = (\varphi_1(x, p_x; T_1), \varphi_2(x, p_x; T_1)),$$

where the maps $\varphi$ and the time $T_1$ are as before, while $H_2$ is defined by

$$H_2(x, p_x) = (\tilde{\varphi}_1(x, p_x; T_1), \tilde{\varphi}_2(x, p_x; T_1)),$$

where $\tilde{\varphi}(x, p_x; \cdot) : \mathbb{R} \to \mathbb{R}^4$ is the solution of the equation (1.1) with initial conditions $x(0) = x$, $\dot{x}(0) = p_x$, $y(0) = 0$, $\dot{y}(0) = -\sqrt{2\Omega(x, 0) + h - p_x^2}$ and $T_1 > 0$ is such that $\varphi_3(x, p_x; T_1) = 0$ and $\varphi_3(x, p_x; t) < 0$ for all $t \in (0, T_1)$. Of course, $D(P) \subset D(H_1)$, $H_1(D(P)) \subset D(H_2)$ and $P = H_2 \circ H_1$. In the following we refer to the maps $H_1$ and $H_2$ as the first and second half Poincaré maps respectively. In other words the map $H_1$ (resp. $H_2$) is the transition map from the section $y = 0$, $p_y > 0$ to the section $y = 0$, $p_y < 0$ (resp. from the section $y = 0$, $p_y < 0$ to the section $y = 0$, $p_y > 0$).

By inspection it is easy to see that, if $(x(t), y(t))$ is a solution of (1.1), then $(\tilde{x}(t), \tilde{y}(t)) := (x(-t), -y(-t))$ is also a solution. This implies that $H_1(x_1, p_1) = (x_2, p_2)$ is equivalent to $H_2(x_2, -p_2) = (x_1, -p_1)$ and $P(x_1, p_1) = (x_2, p_2)$ is equivalent to $P(x_2, -p_2) = (x_1, -p_1)$.

## 4. Periodic orbits

A standard method for studying periodic orbits consists in looking for fixed points of the Poincaré map $P$. By the symmetry of the system considered at the end of the previous section we infer that $H_1(x_1, 0) = (x_2, 0)$ yields $H_2(x_2, 0) = (x_1, 0)$, which in turn implies that $(x_1, 0)$ is a fixed point for the Poincaré map. On the other hand, if the system admits a periodic orbit which crosses orthogonally the $x$-axis at some point $x_1$, then, by definition of Poincaré map, $P(x_1, 0) = (x_1, 0)$, and this is possible only if $H_1(x_1, 0) = (x_2, 0)$ for some $x_2$. It turns out that it is also useful to consider whether there exist points $x_1$ and $x_2$ such that $P(x_1, 0) = (x_2, 0)$, with $x_1 \neq x_2$. By the same reason as before, these points correspond to periodic points of period 2 for the Poincaré map, hence to periodic trajectories for the system crossing the $y = 0$ hyperplane at two different points in each direction, see the pictures of the orbits $D_1$ and $D_2$ below.

Let $f(x)$ be the second component of $H_1(x, 0)$ and $g(x)$ be the second component of $P(x, 0)$. In order to find fixed or period 2 points for the Poincaré map we can look for zeros of the function $f$ or $g$. We remark that, even by considering the map $g$ only, one can still find all the fixed points for the Poincaré map. On the other hand it may happen that $P(x_1, 0) = (x_2, 0)$ with $x_1$ very close to $x_2$, in the sense that $|x_1 - x_2|$ is smaller than the numerical error. In this case it is impossible to find out whether $x_1$ is a fixed point or a periodic point of period 2, without considering either the derivative of the Poincaré map or the half Poincaré map. For this reason it is convenient to study both the map $f$ and the map $g$: first we look for zeros of $f$, i.e. fixed points of $P$, then we can look for zeros of $g$ which are not zeros for $f$, such points yielding periodic points of minimal period 2 of $P$. The following picture displays a numerical computation of $f(x)$ with $x \in (-.62, .62)$ and $h = .1$:

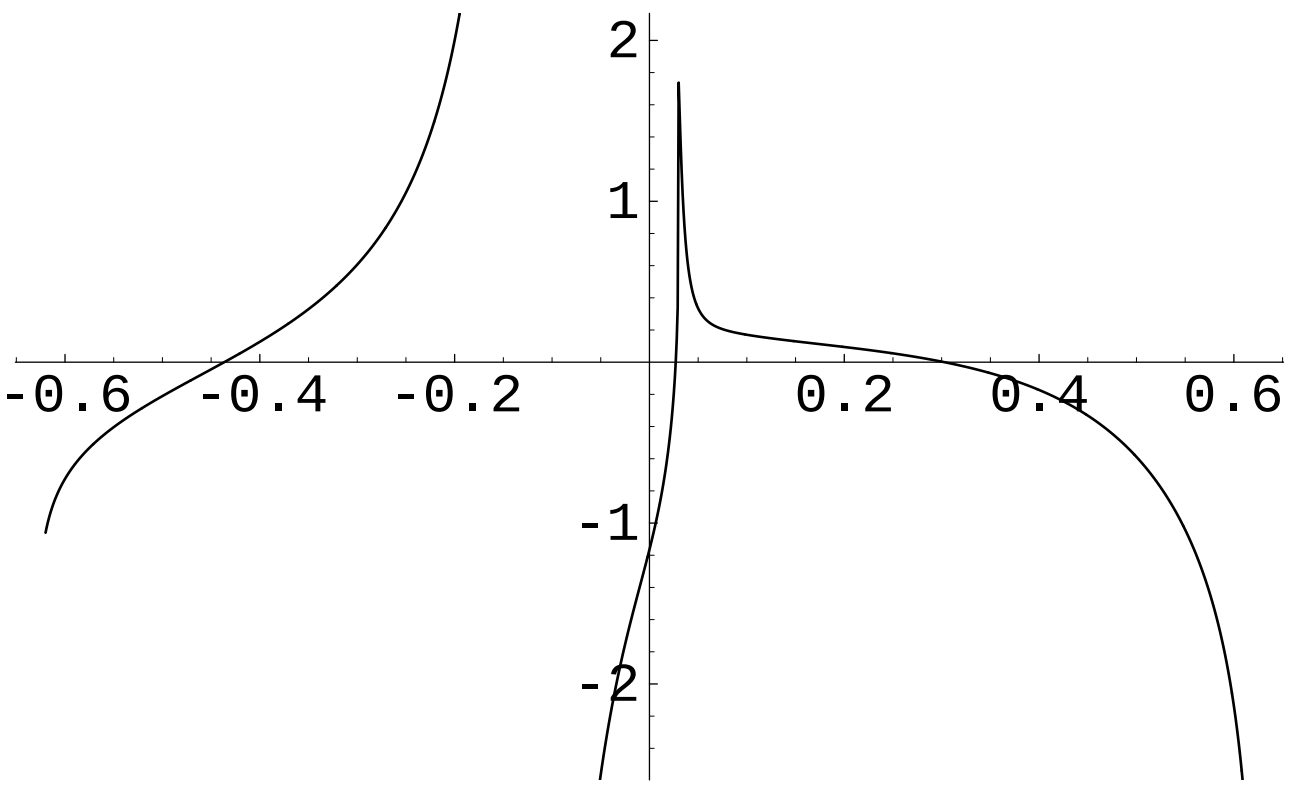

$h = .1$

The picture strongly suggests that the Poincaré map has three fixed points in the set $\{(x, 0) : x \in (-.62, .62)\}$, which is exactly the result we show rigorously in the next section, where we also give a narrow bound on the position of such points. Our strategy to detect and prove the existence of periodic orbits is as follows: first we choose some values for $h$ and compute an approximate image of the $x$ axis through the maps $f$ and $g$ as in the previous picture; in this way we can spot the places where the intersections should be located. Assume that the numerically computed graph of $f$ crosses the $x-$axis in a neighborhood of some point $\bar{x}$: we conjecture the existence of a fixed point for $P$ nearby some point $(\bar{x}, 0)$. To prove the conjecture we choose $x_1$ and $x_2$ such that $x_1 < \bar{x} < x_2$ and both $x_1$ and $x_2$ are very close to $\bar{x}$. Then we compute the rigorous half Poincaré map $H_1$ at $(x_1, 0)$ and $(x_2, 0)$. If we can prove that the second component of $H_1(x_1, 0)$ has opposite sign with respect to the second component $H_1(x_2, 0)$ and that the segment joining the two points belongs entirely to the domain of $H_1$, then by the continuity of the half Poincaré map we have proved that there exists at least a point $x_1 < \tilde{x} < x_2$ such that $H_1(\tilde{x}, 0)$ lies on the $x$ axis, therefore a periodic orbit intersects the $x$-axis orthogonally at $(\tilde{x}, 0)$. On the other hand, if we prove that some portion of the $x$-axis is mapped away from the $x$-axis itself, then we have a proof that there are no periodic solutions which cross the $x$-axis orthogonally in that section. Then we search for points of period 2 for the map $P$, that is we study the numerically computed graph of $g$ and look for intersections with the $x-$axis. If we find two points $x_1$, $x_2$ such that the second components of $P(x_1, 0)$ and $P(x_2, 0)$ lie on the opposite sides of the $x$-axis, the set $[x_1, x_2] \times \{0\}$ belongs to the domain of $P$ and $P([x_1, x_2], 0)$ does not intersect the set $[x_1, x_2] \times \{0\}$, then we have a proof that there exists at least a point $x_1 < \tilde{x} < x_2$ such that $\tilde{x}$ is the intersection of an orbit of minimal period 2 with the $x$-axis.

The case we consider is usually referred to as "the Copenhagen orbits", from the results of the Observatory of Copenhagen (see [S] and references therein). We recall that, by our choice of the constant $C$ in (2.2), $h = 0$ is the energy of the stationary solution at the Lagrangian point $L_1$ (the origin). Since the problem has the additional symmetry consisting in switching the primaries, we only look for orbits that cross the portion of the $x$-axis between the primaries with positive speed in the $y$ direction. The lowest value of $h$ we consider is 0, when the admissible region is split in two parts touching at the origin and no trajectory can enter both regions. The highest value of $h$ we consider is 0.8, since at slightly larger value the bounded admissible region touches the unbounded part and some trajectory starting close to the primaries may be unbounded.

| $h$ | $S_1$ | $S_2$ | $L$ | $D_1$ | $D_2$ | $U_1$ | $U_2$ |
|---|---|---|---|---|---|---|---|
| 0 | 0.3158 | −0.4399 | | | | | |
| .1 | 0.3028 | −0.4365 | 0.02697 | | | | |
| .2 | 0.2883 | −0.4330 | 0.03894 | | | | |
| .22 | 0.2851 | −0.4323 | 0.04102 | | | | |
| .24 | 0.2818 | −0.4316 | 0.04303 | | | 0.06963 | 0.1043 |
| .26 | 0.2784 | −0.4309 | 0.04497 | 0.04423 | 0.04456 | 0.06737 | 0.1205 |
| .28 | 0.2749 | −0.4301 | 0.04687 | 0.04600 | 0.04651 | 0.06689 | 0.1343 |
| .3 | 0.2712 | −0.4294 | 0.04873 | 0.04775 | 0.0484 | 0.06712 | 0.1470 |
| .4 | 0.2443 | −0.4257 | 0.05751 | 0.05615 | 0.05730 | 0.07199 | 0.2085 |
| .5 | | −0.4218 | 0.06575 | 0.06412 | 0.06550 | 0.07870 | |
| .6 | | −0.4178 | 0.07369 | 0.07185 | 0.07345 | 0.08595 | |
| .8 | | −0.4093 | 0.08922 | 0.08716 | | | |

Table 1

The main result of this section is that the system (1.1) admits periodic solutions as shown in Table 1. In the left hand side column the energy level is displayed, while the remaining columns represent the $x$-coordinate $\pm 5 \cdot 10^{-4}$ of the intersection of the orbits with the portion of the $x$-axis between the primaries (the $L$, $D_1$ and $D_2$ orbits have two such intersections; for the $L$ orbits we only consider the intersection with positive $y$ velocity, the other being symmetric with respect to the origin; for the $D_1$ and $D_2$ orbits we consider the only orthogonal intersection). More precisely, the following theorem holds:

**Theorem 4.1.** *For all energy values displayed in column $h$ of Table 1, the system (1.1) admits at least a periodic solution for each value printed in the remaining columns. Such solution crosses the $x$-axis orthogonally twice. The $x$-coordinate of the intersection with positive $y$-velocity lies in the interval centered in the position given in the table with width $10^{-3}$. The solutions corresponding to values in the columns $S_1$, $S_2$, $L$, $U_1$ and $U_2$ do not intersect the $x$-axis at any other point, while the solutions corresponding to values in the columns $D_1$ and $D_2$ intersect the $x$-axis twice at another point with the same negative $x$-velocity and with opposite $y$-velocity. The orbits in the column $S_2$ are retrograde around $P_1$; the orbits in the classes $S_1$, $U_1$ and $U_2$ are direct around $P_2$. The orbits in the class $L$ are retrograde around $L_1$.*

*Proof. Orbits $S_1$, $S_2$, $L$, $U_1$ and $U_2$.* We checked by computer assistance (see Section 7) that the second components of $H_1(x + 5 \cdot 10^{-4}, 0)$ and $H_1(x - 5 \cdot 10^{-4}, 0)$ have different sign and the interval set $[x - 5 \cdot 10^{-4}, x + 5 \cdot 10^{-4}] \times \{0\}$ is in the domain of $H_1$, where $x$ is any of the values in the columns $S_1$, $S_2$, $L$, $U_1$ and $U_2$. According to the argument presented at the beginning of this section, this suffices to prove the existence of a fixed point of the Poincaré map $P$ and hence the existence of a periodic orbit. To prove that an orbit in the column $S_2$ is retrograde around $P_1$ we proceed as follows. We compute the trajectory of the set $[x - 5 \cdot 10^{-4}, x + 5 \cdot 10^{-4}] \times \{0\}$ until it crosses the Poincaré section and we check that the angular velocity with respect to $P_1$ is always strictly positive. This also implies that the $x$-axis is crossed only once in each direction, and by construction and the symmetry of the system the intersections are othogonal. To prove that the orbits in the classes $S_1$, $U_1$ and $U_2$ are direct around $P_2$, the orbits in the class $L$ are retrograde around $L_1$ and they also cross the $x$-axis only once in each direction we follow an analogous procedure.

*Orbits $D_1$ and $D_2$.* Same argument with the map $P$. To check that the periodic points of the Poincaré map corresponding to these orbits have minimal period 2 we checked that $P([x - 5 \cdot 10^{-4}, x + 5 \cdot 10^{-4}], 0)$ is mapped away from $[x - 5 \cdot 10^{-4}, x + 5 \cdot 10^{-4}] \times \{0\}$. The symmetry of the system forces the two non-orthogonal intersections with the $x$-axis to occur at the same point, with the same $x$-velocity and with opposite $y$-velocity. □

**Remark 4.1.** *Orbits belonging to the same column appear to belong to the same class, in the sense that they have the same linear stability and they have the same winding number with respect to the primaries and the Lagrangian point $L_1$. More precisely, all orbits are unstable, except for orbits $S_1$ and $S_2$ which are linearly stable. The result on the stability is purely numerical.*

The orbits in the class $L$ are well known, their trajectory is very close to an ellipse with center at the Lagrangian point $L_1$. They branch out from $L_1$, in the sense that as $h \to 0$ they collapse to $L_1$. These are the only orbits considered in this paper which do not wind around any primary, while they wind around $L_1$. The orbits in the class $D_1$ and $D_2$ correspond to points of period 2 of the Poincaré map. In the following pictures the orbits $D_1$, $D_2$ and $L$ at energy level $h = .5$ are displayed (the small disk represents $P_1$).

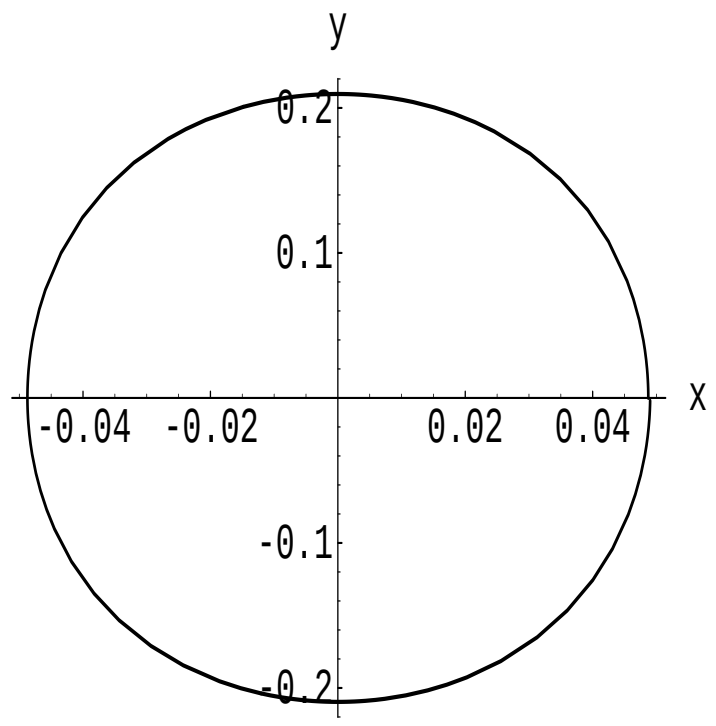

Orbit $L$

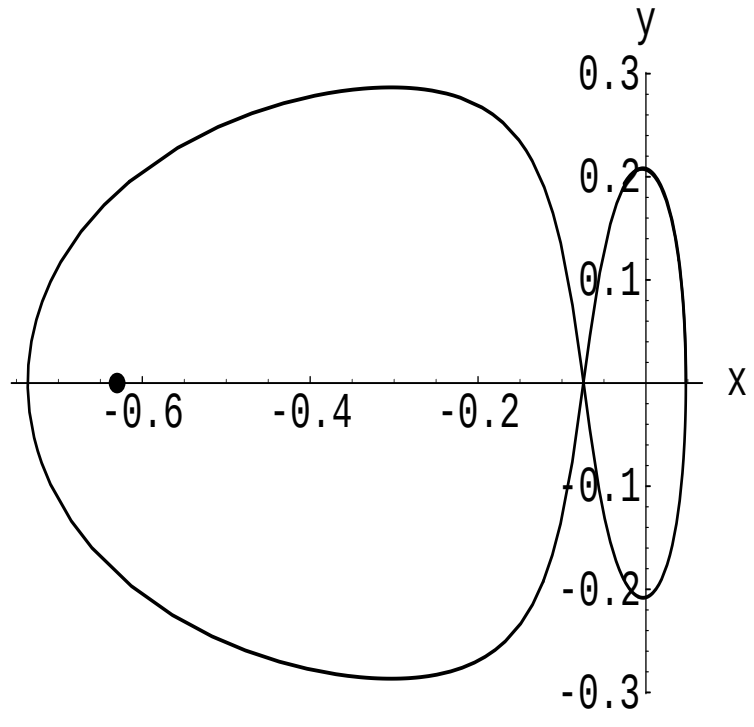

Orbit $D_1$

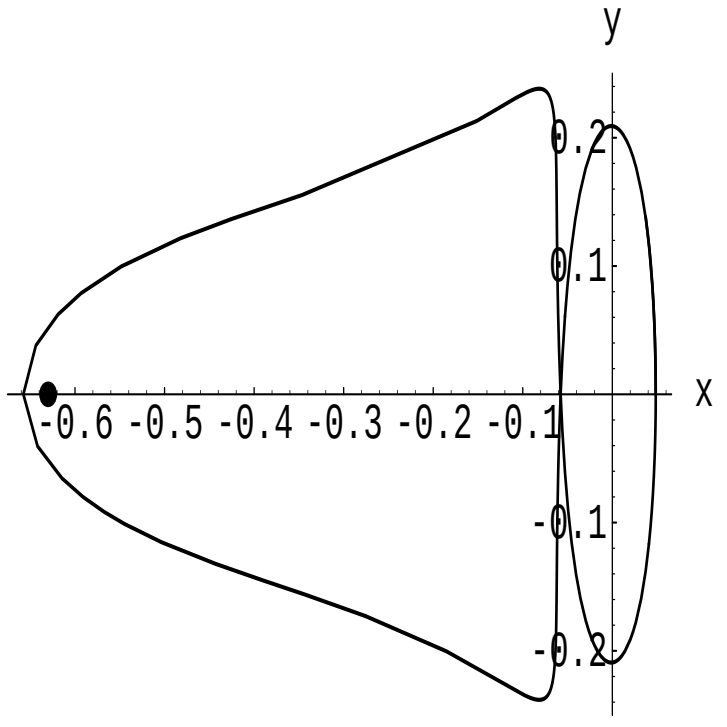

Orbit $D_2$

The same orbits are plotted together in Figure 1: note that part of the trajectories of both $D_1$ and $D_2$ are very close to the trajectory of $L$, giving a visual image of the strong instability of the system.

**Remark 4.2.** *The existence of unstable orbits around the other Lagrangian points is also well-known, but they occur only for higher values of the energy; we do not investigate those orbits, but we claim that the methods presented in this paper could be used to study those orbits and, possibly, chaotic dynamics involving those orbits as well.*

In Figure 2 we plot all periodic orbits at energy level $h = .3$: note that the shape of the orbits $D_1$, $D_2$ and $L$ are quite similar to the same orbits in Figure 1, supporting the claim that they belong to the same branch.

## 5. Chaotic dynamics

We apply the method developed in [Z1, Z2, Z3], see also [AZ] where the Hénon-Heiles Hamiltonian was concerned. We only consider the system at energy $h = .3$. We prove the existence of symbolic dynamics on 14 symbols and we give a lower estimate for the topological entropy.

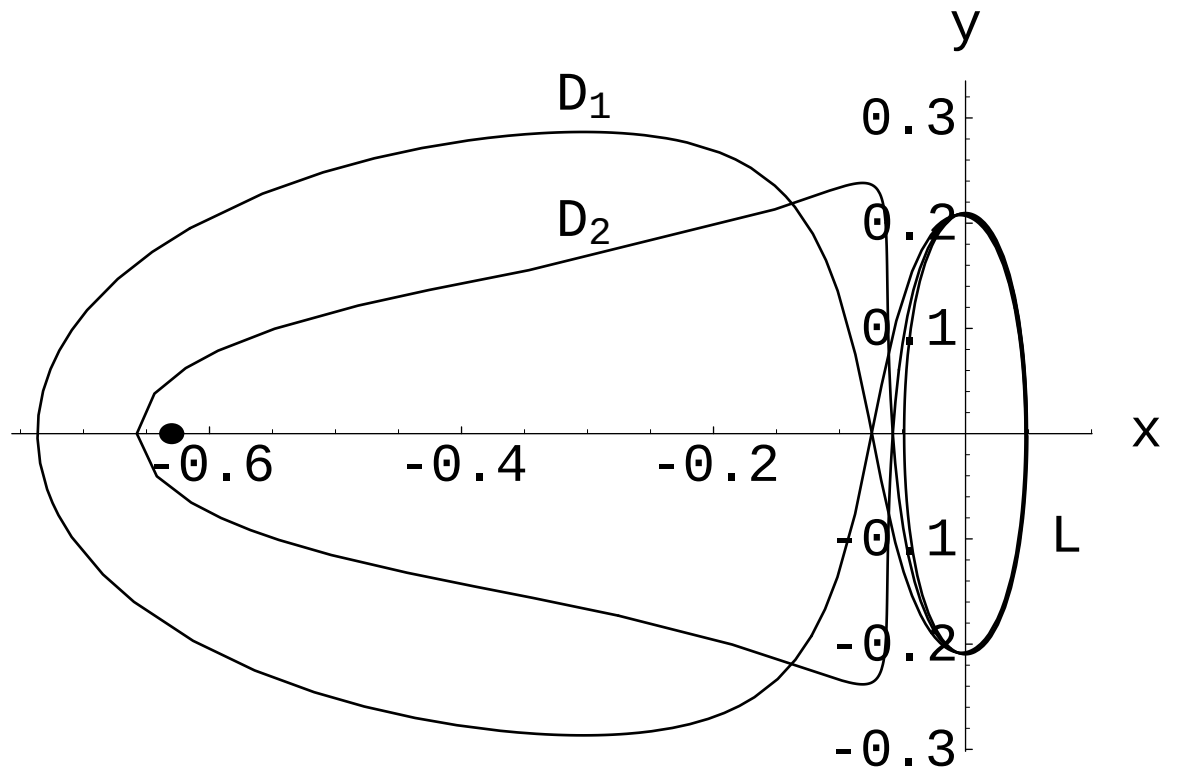

FIGURE 1. Orbits $D_1$, $D_2$ and $L$ at energy level $h = .5$

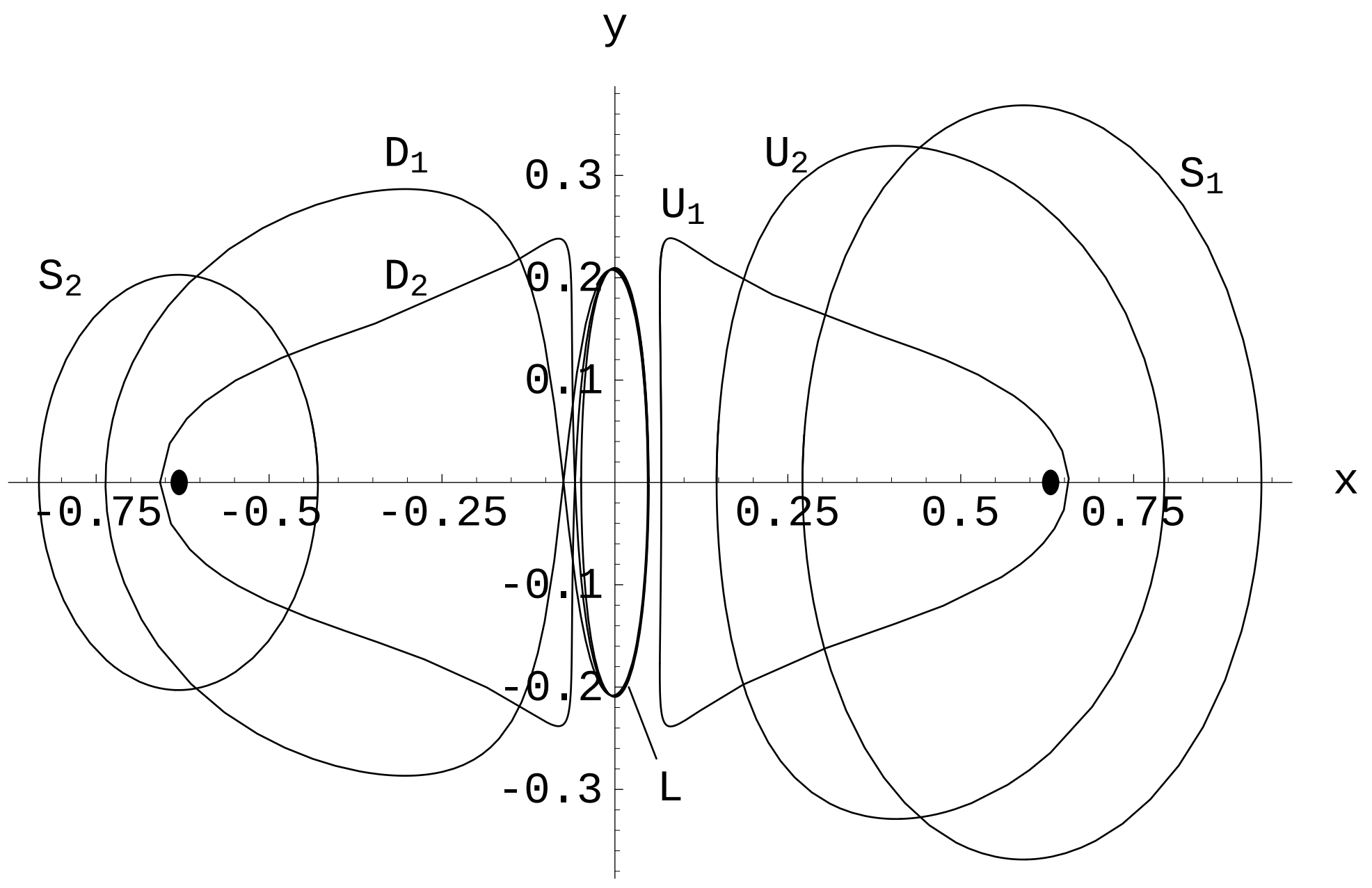

FIGURE 2. All orbits at $h = .3$

The computation of rigorous estimates in the (restricted) 3-body problem presents more difficulties than the cases previously considered. The first is due to the fact that the Poincaré map is not defined in a connected region, indeed we have to exclude at least the lines $x = R_1$ and $x = -R_2$. This could be avoided by using some kind of regularization, but we prefer to keep the original coordinates, both because they are far more intuitive and because the Levi-Civita transform (or its variations) is not one to one. Furthermore, even after some kind of regularization, it is far from trivial to determine what is the domain of the Poincaré map. But the main problem is due to the fact that, even using a sofisticated algorithm to compute the rigorous bounds, such bounds turn out to be very large, particularly in the $p_x$ direction. The reason for such large bounds are not only the inevitable computational errors, which in theory can be as small as we like (by taking a smaller time step at the price of a slower computation), but particularly the wrapping effect. See [AZ, GZ] and the references given there for a discussion on this topic. Here we just want to point out that although it is possible by a change of variable to get rid of the singularities, it is not possible to get rid of the wrapping effect in the same way, at least not with the same change of variable.

The first trick we had to adopt consists in computing only half of the Poincaré map at a time, and this is the second reason, in fact the most important, for introducing the maps $H_1$ and $H_2$. Indeed the wrapping effect is usually exponential, therefore by considering about half trajectory time it is rather drastically reduced. By this method we obtain a great reduction of the error: of course we have to pay the price of a larger number of computations, but the trade-off is very positive. The second trick that proved to be essential consists in computing the inverse of the map instead of the actual map for some checks, see the definition of back-covering below. This is important whenever some trajectory starts away from the primaries but comes close to one of them at the intersection with the Poincaré plane, indeed in such cases a very short time step is necessary to keep the error small, since close to the singularities the speed undergoes a strong variation while the particle crosses the Poincaré plane. The particle is in fact a point, but in order to compute the Poincaré map we have to consider the envelope of its position with the error bounds, hence it takes a finite time to cross the plane. It turns out that the inverse trajectory, starting close to a primary and ending away from both primaries, raises a much lower error.

5.1. **Topological tools.** Definitions 5.1 and 5.3 were introduced in [AZ]. Definitions 5.4 and 5.5 are introduced here for the first time to deal with the new difficulties of the restricted 3-body problem. Lemma 5.3 and Corollaries 5.4 and 5.5 are also new.

**Definition 5.1.** *A **triple set** (or t-set) is a triple $N = (|N|, N^l, N^r)$ of closed subsets of $\mathbb{R}^2$ satisfying the following properties:*

- **1:** *$|N|$ is a closed parallelogram*
- **2:** *$N^l$ and $N^r$ are closed half-planes*
- **3:** *the sets $N^{le} := N^l \cap |N|$ and $N^{re} := N^r \cap |N|$ are two nonadjacent edges of $|N|$*

*We call $|N|$, $N^l$, $N^r$, $N^{le}$ and $N^{re}$ the support, the left side, the right side, the left edge and the right edge of the t-set $N$ respectively. One observes that $\mathbb{R}^2 \setminus (|N| \cup N^l \cup N^r)$ consists of two disjoint sets. We call $N^t$ and $N^b$ the closure of such sets (the top and bottom sides of the triple set). We also set $N^{te} := N^t \cap |N|$ and $N^{be} := N^b \cap |N|$ (the top and bottom edges of $N$).*

**Remark 5.1.** Since all theorems we use concern topological properties of the t-sets, then one can choose any continuous deformation of the t-set defined above, obtaining the same results. An example of the sets we actually use, whose exact definition is given in Section 7, is given in Figure 3.

To exploit the symmetry of the system we give the following definition:

**Definition 5.2.** *Let $M$ be a t-set. We define its symmetric image with respect to the $x$-axis $\tilde{M}$ as follows: if $S : \mathbb{R}^2 \to \mathbb{R}^2$ is the map defined by $S(x, y) = (x, -y)$, let $|\tilde{M}| = S(|M|)$, $\tilde{M}^l = S(M^t)$, $\tilde{M}^r = S(M^b)$, and the remaining definitions follow as in Definition 5.1.*

**Definition 5.3.** *Let $f : \Omega \subset \mathbb{R}^2 \to \mathbb{R}^2$ be a map and let $N_1$ and $N_2$ be two triple sets. We say that $N_1$ $f-$covers $N_2$ ($N_1 \overset{f}{\Longrightarrow} N_2$) if:*

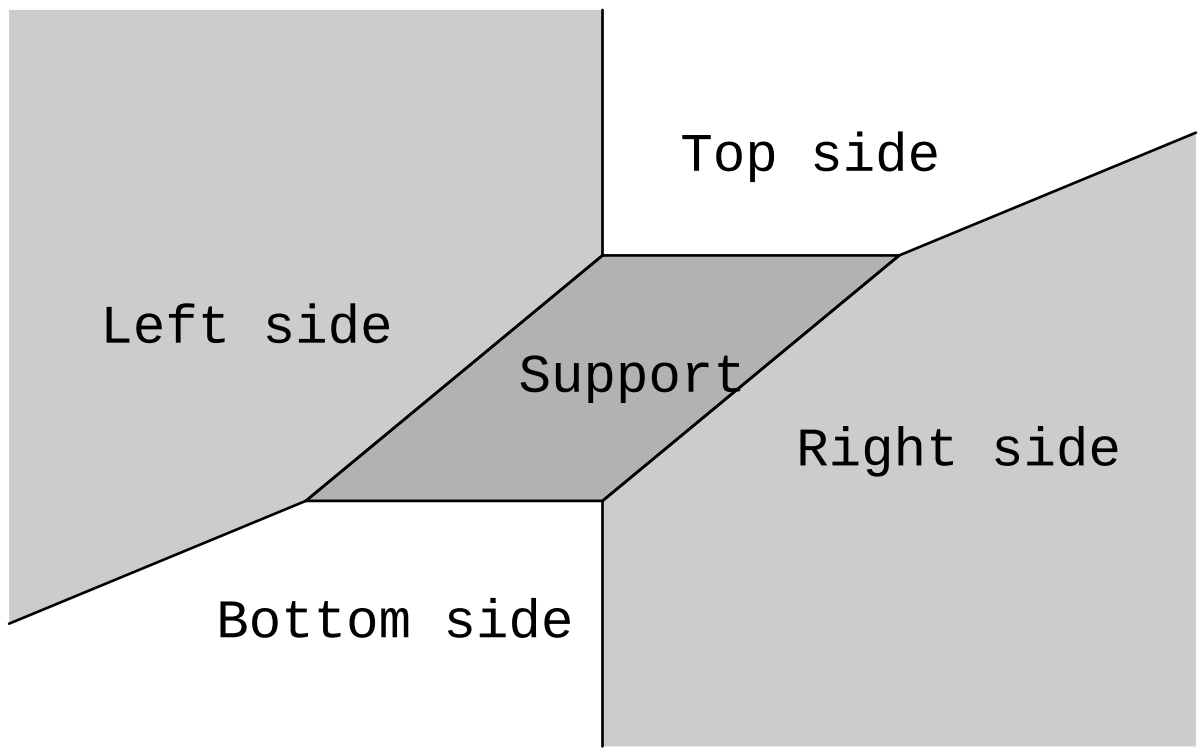

FIGURE 3. An example of a t-set

**a:** $f(|N_1|) \subset \mathrm{int}(N_2^l \cup |N_2| \cup N_2^r)$
**b:** *either* $f(N_1^{le}) \subset \mathrm{int}(N_2^l)$ *and* $f(N_1^{re}) \subset \mathrm{int}(N_2^r)$
*or* $f(N_1^{le}) \subset \mathrm{int}(N_2^r)$ *and* $f(N_1^{re}) \subset \mathrm{int}(N_2^l)$

The following lemma says that we can reduce the condition (a) in the above definition to the boundary of $|N_1|$ if we know that the map $f$ is defined on $|N_1|$ and it is injective. This lemma plays a very important role in the computer assisted verification of the covering relations, as it reduces the computations to the boundary of $|N_1|$ (see Section 6 in [GZ], for more details).

**Lemma 5.1.** *Let* $f : \Omega \subset \mathbb{R}^2 \to \mathbb{R}^2$ *be a map and let* $N_1$ *and* $N_2$ *be two triple sets. Assume that* $f$ *is an injective map on* $|N_1|$*, then* $N_1 \overset{f}{\Longrightarrow} N_2$ *if and only if*
**a′:** $f(\partial|N_1|) \subset \mathrm{int}(N_2^l \cup |N_2| \cup N_2^r)$
**b:** *either* $f(N_1^{le}) \subset \mathrm{int}(N_2^l)$ *and* $f(N_1^{re}) \subset \mathrm{int}(N_2^r)$
*or* $f(N_1^{le}) \subset \mathrm{int}(N_2^r)$ *and* $f(N_1^{re}) \subset \mathrm{int}(N_2^l)$

As we pointed out above, in some circumstances it is easier to compute the inverse flow, then the direct flow. Furthermore, we need to deal with the half Poincaré maps. To this purpose we give the following definitions:

**Definition 5.4.** *Let* $N_1$ *and* $N_2$ *be two triple sets. We say that* $N_1$ $f-$*backcovers* $N_2$ *(*$N_1 \overset{f}{\Longleftarrow} N_2$*) whenever:*
**a:** $f : \Omega_1 \subset \mathbb{R}^2 \to \Omega_2 \subset \mathbb{R}^2$ *is a homeomorphism*
**b:** $|N_2| \subset \Omega_2$
**c:** $f^{-1}(\partial|N_2|) \subset \mathrm{int}(N_1^t \cup |N_1| \cup N_1^b)$
**d:** *either* $f^{-1}(N_2^{te}) \subset \mathrm{int}(N_1^t)$ *and* $f^{-1}(N_2^{be}) \subset \mathrm{int}(N_1^b)$
*or* $f^{-1}(N_2^{te}) \subset \mathrm{int}(N_1^b)$ *and* $f^{-1}(N_2^{be}) \subset \mathrm{int}(N_1^t)$

**Definition 5.5.** *Let* $N_1$ *and* $N_2$ *be two triple sets. We say that* $N_1$ *generically* $f-$*covers* $N_2$ *(*$N_1 \overset{f}{\Longleftrightarrow} N_2$*) if* $N_1$ $f-$*covers* $N_2$ *or there exists* $n \geq 1$ *t-sets* $M_i$*,* $i = 1, \ldots, n$ *and* $n+1$ *maps* $g_j$*,* $j = 0, \ldots, n$ *such that* $f = g_j \circ \cdots \circ g_1 \circ g_0$*,* $N_1 \overset{g_0}{\Longrightarrow} M_1$*,* $M_i \overset{g_i}{\Longrightarrow} M_{i+1}$ *for all* $i = 1, \ldots, n-1$ *and* $M_n \overset{g_n}{\Longrightarrow} N_2$*.*

**Remark 5.2.** *Although covering, backcovering and generic covering occur often simultaneously and they are indeed a very similar phenomenon, they are not equivalent. In fact it may even happen that a map $f$ is not defined on the whole support of $N_1$, and still $N_1$ $f-$backcovers $N_2$ or $N_1$ generically $f-$covers $N_2$ . On the other hand, the result of a backcovering or of a generic covering relation is very similar to the result of a covering relation as far as the results we are interested in are concerned, see Corollaries 5.4 and 5.5.*

The following theorem follows immediately from Theorem 4 in [Z3] and together with Theorem 5.6 it is the main topological tool we use:

**Theorem 5.2.** *Given $n$ t-sets $M_i \subset \mathbb{R}^2$ and $n$ continuous maps $f_i : M_i \to \mathbb{R}^2$, such that*

$$M_0 \stackrel{f_0}{\Longrightarrow} M_1 \stackrel{f_1}{\Longrightarrow} M_2 \stackrel{f_2}{\Longrightarrow} M_2 \ldots \stackrel{f_{n-1}}{\Longrightarrow} M_0 = M_n,$$

*then there exists $x \in \operatorname{int}|M_0|$, such that $f_k \circ \cdots \circ f_1 \circ f_0(x) \in \operatorname{int}|M_{k+1}|$, for $k = 0, \ldots, n-1$ and $x = f_{n-1} \circ \cdots \circ f_1 \circ f_0(x)$.*

By the following lemma we can extend the previous theorem to the generic backcovering, see Corollary 5.5.

**Lemma 5.3.** *Let $\Omega_1$ and $\Omega_2$ be two open sets of $\mathbb{R}^2$, let $f : \Omega_1 \to \Omega_2$ be a homeomorphism and let $N_1$ and $N_2$ be two triple sets such that $|N_i| \subset \Omega_i$, $i = 1, 2$. If $N_1$ $f-$backcovers $N_2$ then there exists a t-set $K \subset \Omega_1$ such that $N_1$ id$-$covers $K$, and $K$ $f-$covers $N_2$ (id is the identity map in $\Omega_1$).*

*Proof.* By definition of backcovering, if $N_1$ $f-$backcovers $N_2$, then $f^{-1}(|N_2|)$ looks as the light grey rectangle in the Figure 4, therefore it is possible to define a t-set $K$ such that the boundary of its support is very close to the boundary of $f^{-1}(|N_2|)$, in such a way that $N_1 id-$covers $K$ and $K$ $f-$covers $N_2$. To prove this, just consider a set $K'$ as the dark grey rectangle in the picture or, more precisely, if $a$ and $b$ are the length of the sides of $|N_2|$, let $K'$ be the parallelogram with the same center as $|N_2|$, the sides parallel to the sides of $|N_2|$ and such that the length of its sides is $(1-\varepsilon)a$ and $(1+\varepsilon)b$, $0 < \varepsilon < 1$. If $\varepsilon$ is small enough, $K' \subset \Omega_2$ and $f^{-1}(K')$ is defined. Call the edges of $K'$ top, bottom, left and right according to the corresponding edge of $N_2$. Now define $K$ as follows. The edges are the counterimage of the edges of $K'$. By construction, definition of backcovering, compactness of $K'$ and continuity of $f^{-1}$, if $\varepsilon$ is small enough the top (resp. bottom) edge of $K$ lies in the top (resp. bottom) side of $N_1$ and both the right and the left edges are in $\operatorname{int}(N_1^t \cup |N_1| \cup N_1^b)$, hence it is possible to define the left (resp. right) side of $K$ in such a way that $N_1$ $id-$covers $K$. By construction $K$ $f-$covers $N_2$ and the proof is complete. □

**Corollary 5.4.** *Let $N_1$ and $N_2$ be two triple sets and let $f$ be a homeomorphism. If $N_1$ $f-$backcovers $N_2$, then $N_1$ generically $f-$covers $N_2$.*

**Corollary 5.5.** *Given t-sets $M_i \subset \mathbb{R}^2$ and $n$ continuous maps $f_i : \Omega_{1i} \subset \mathbb{R}^2 \to \Omega_{2i} \subset \mathbb{R}^2$ such that*

$$M_0 \stackrel{f_0}{\Longleftrightarrow} M_1 \stackrel{f_1}{\Longleftrightarrow} M_2 \stackrel{f_2}{\Longleftrightarrow} M_2 \ldots \stackrel{f_{n-1}}{\Longleftrightarrow} M_0 = M_n,$$

*then there exists $x \in \operatorname{int}|M_0|$, such that $f_k \circ \cdots \circ f_1 \circ f_0(x) \in \operatorname{int}|M_{k+1}|$, for $k = 0, \ldots, n-1$ and $x = f_{n-1} \circ \cdots \circ f_1 \circ f_0(x)$.*

Assume that we have $n$ t-sets $N_i$, $i = 1, \ldots, n$, with some covering relations. Let $N = \bigcup_i |N_i|$; the following definitions are standard.

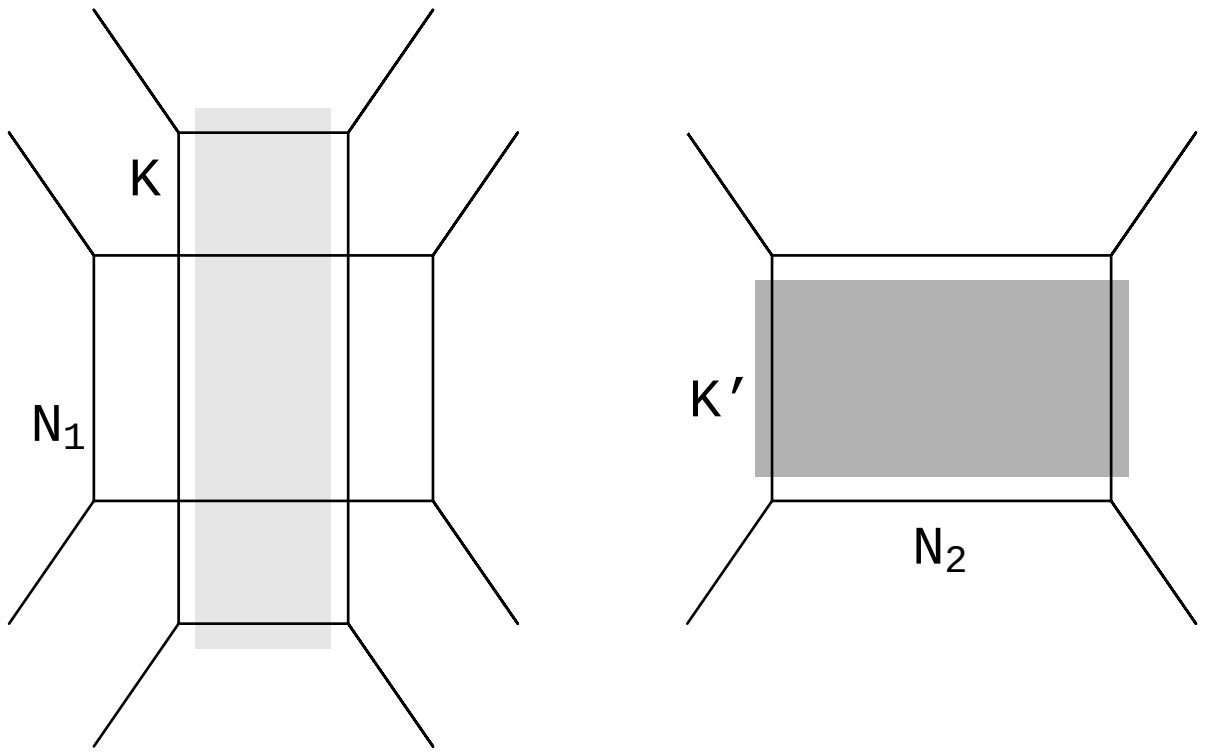

FIGURE 4. In light grey $f^{-1}(N_2)$, in dark grey $K' = f(K)$

**Definition 5.6.** *Let $f$ be injective. The invariant set of $N$ is defined by* $\mathrm{Inv}(N, f) := \{x \in N : f^i(x) \in N$ *for all* $i \in \mathbb{Z}\}$.

**Definition 5.7.** *The transition matrix $T(j,i)$, $i, j = 1, \ldots, n$, is defined as follows:*

$$T(j,i) = \begin{cases} 1 & \text{if } N_i \overset{f}{\Longleftrightarrow} N_j \\ 0 & \text{otherwise} \end{cases}$$

Let $\Sigma_n$ be the set of bi-infinite sequences of $n$ symbols

**Definition 5.8.** *A sequence $\{x_k\} \in \Sigma_n$ is said to be* admissible *if $T(x_{k+1}, x_k) = 1$ for all $k$. We denote by $\Sigma_A \subset \Sigma_n$ the set of all admissible sequences.*

**Definition 5.9.** *Assume $|N_i| \cap |N_j| = \emptyset$, for $i \neq j$. The projection* $\pi : \mathrm{Inv}(N, f) \to \Sigma_A$ *is defined by setting $\pi(x)_i = j$ where $j$ satisfies $f^i(x) \in |N_j|$ for all $i \in \mathbb{Z}$.*

The set $\Sigma_A$ inherits the topology from $\Sigma_n$; the shift map $\sigma : \Sigma_A \to \Sigma_A$ is continuous. We prove a semiconjugacy between $\sigma$ and $f$, i.e. we prove that $\sigma \circ \pi = \pi \circ f|_{\mathrm{Inv}(N,f)}$. In particular this implies that there exists a symbolic dynamics structure on $\mathrm{Inv}(N, f)$.

The following theorem was proved in [Z3] (see Theorems 5 and 6) for the case $n = 2$. The following is a natural extension to a generic number of sets and the proof is exactly the same.

**Theorem 5.6.** *The projection $\pi$ is onto, and if $\{x_n\} \in \Sigma_A$ is a periodic sequence, then $\pi^{-1}(\{x_n\})$ contains a periodic point.*

We point out that these kind of topological tools have been created to deal with hyperbolic periodic points and are not suitable to look for elliptic points. On the other hand the symmetry of this system allows the search and proof of (symmetric) periodic points in a much easier way, as pointed out in Section 4. Furthermore the main purpose of the paper is to show the existence of a symbolic dynamics, which cannot occur in neighborhoods of stable points, therefore we will not address this topic any further.

5.2. **Heuristic results.** In this section of the paper we provide the heuristic results we obtained, while in the following section we prove that such results are rigorous. In the remaining part of this paper we fix $h = .3$.

We need to define 14 t-sets $N_0, \ldots, N_5$, $\tilde{N}_1, \ldots, \tilde{N}_4$, $K_0, \ldots, K_3$: the precise definition is given in Section 7. Here we only point out that the sets $\tilde{N}_1, \ldots, \tilde{N}_4$ are by definition the symmetric image with respect to the $x$-axis of the sets $N_1, \ldots, N_4$ (see Definition 5.2), while the remaining sets are invariant with respect to the same symmetry. In the following figures the supports of the sets $N_0 - N_5$ and their images through the Poincaré map are represented, together with a portion of the images the supports of the sets $K_0$ and $K_2$. We remark that the images of the sets $N_2$, $K_0$ and $K_2$ are very thin and appear on the picture as lines. The sets are displayed in thick lines, while their images are displayed in thin lines We made two pictures for better clarity.

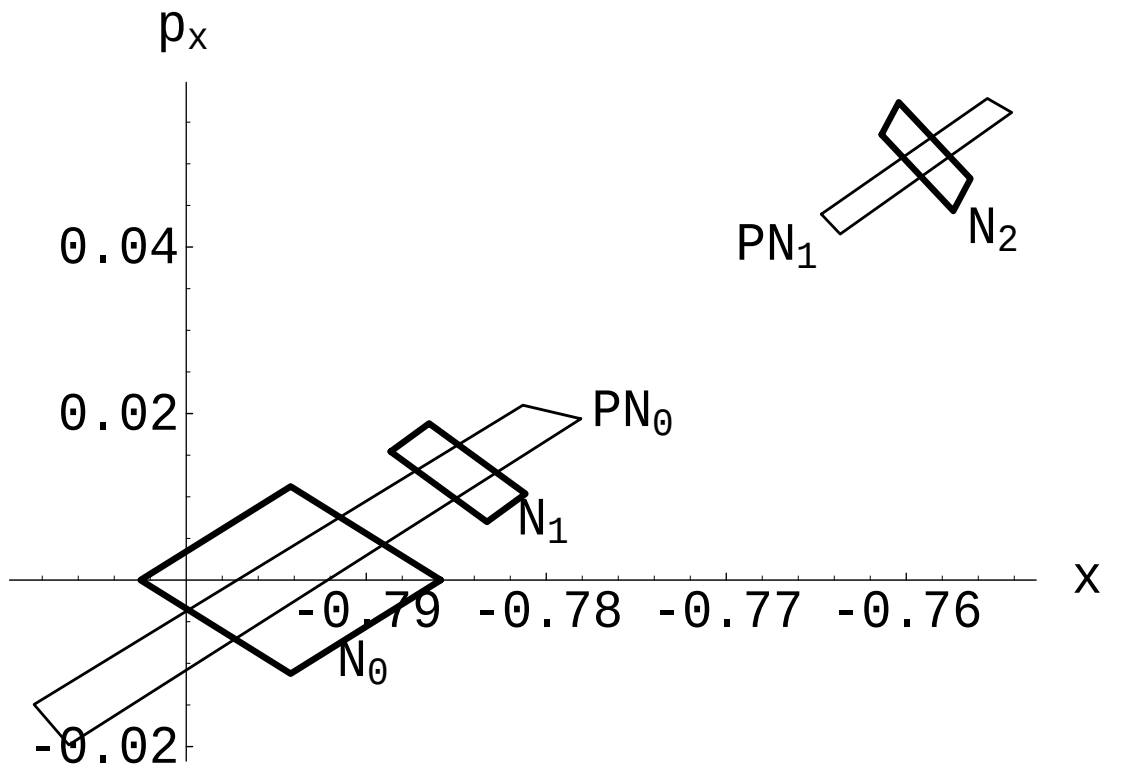

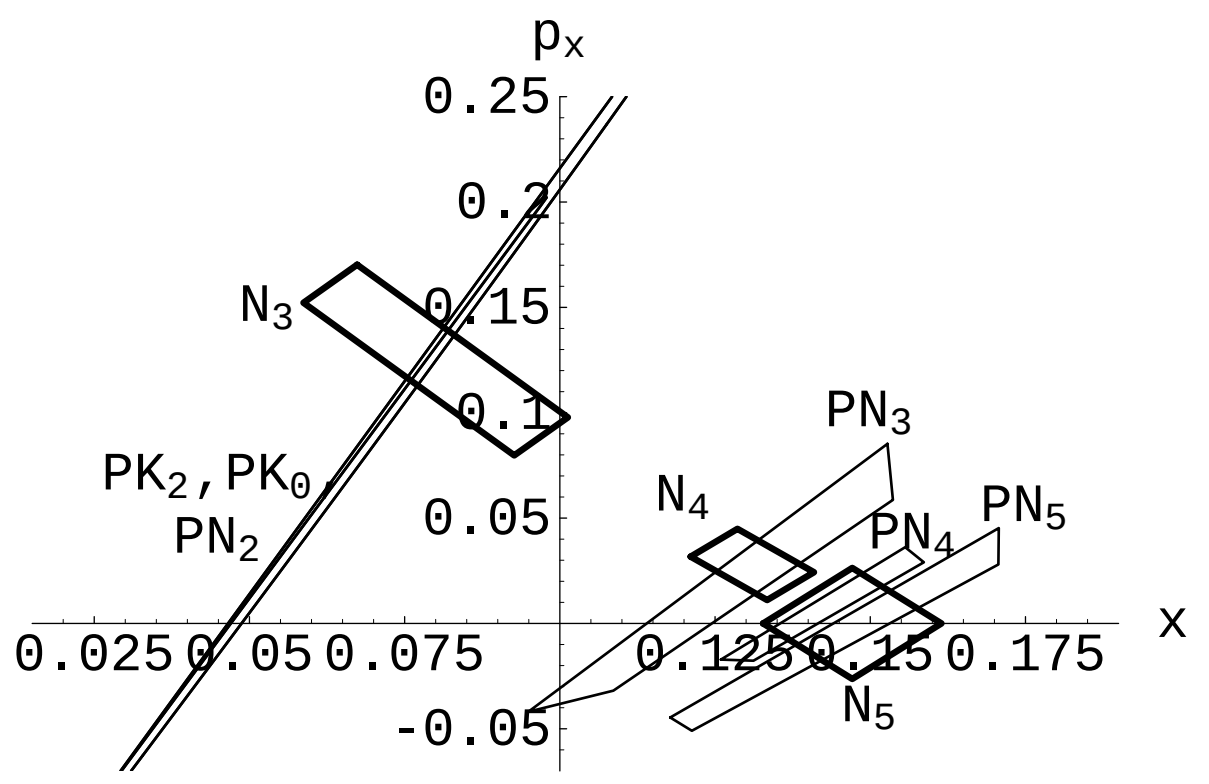

**Lemma 5.7.** *Let $M_1$ and $M_2$ be t-sets, let $\tilde{M}_1$ and $\tilde{M}_2$ be the sets symmetric with respect to symmetry defined above and assume that $M_1 \overset{P}{\Longrightarrow} M_2$. Then $\tilde{M}_2 \overset{P}{\Longleftarrow} \tilde{M}_1$.*

*Proof.* This follows by the Definitions 5.3, 5.4, 5.2 and the symmetry of the Poincaré map (see the end of Section 3). □

In the following we denote by $\Longrightarrow$ (resp. $\Longleftarrow$) the covering (resp. the backcovering) with respect to the map $P$ and we apply the topological theorems introduced above to the map $P$.

The numerical experiments suggest that the following covering relations hold:

$$N_0 \Longrightarrow N_0 \Longrightarrow N_1 \Longrightarrow N_2 \Longrightarrow N_3 \Longrightarrow N_4 \Longrightarrow N_5 \Longrightarrow N_5,$$

$$K_0 \Longrightarrow K_0 \Longrightarrow N_3, \tag{5.1}$$

$$K_1 \Longrightarrow K_2 \Longrightarrow N_3,$$

$$K_3 \Longrightarrow K_3 \Longrightarrow N_3,$$

If these relations could be verified, then by Lemma 5.7 the following covering relations would hold as well:

$$N_5 \Longleftarrow \tilde{N}_4 \Longleftarrow \tilde{N}_3 \Longleftarrow \tilde{N}_2 \Longleftarrow \tilde{N}_1 \Longleftarrow N_0,$$

$$\tilde{N}_3 \Longleftarrow K_0, \tag{5.2}$$

$$\tilde{N}_3 \Longleftarrow K_2 \Longleftarrow K_1,$$

$$\tilde{N}_3 \Longleftarrow K_3.$$

5.3. **Rigorous covering relations.** As we pointed out before, it is not convenient to check directly the covering relations given in the previous section, since that would require an enormous amount of computer time. This problem is due to the fact that, even if we start from a very small interval in the computation of the Poincaré map, the error and the wrapping effect accumulated in the Poincaré time is quite large. Instead we prefer to define 12 auxiliary t-sets $M_i$, $i = 0, \dots, 6$ and $L_i$, $i = 1, \dots, 5$ on the symmetric Poincaré plane, i.e. the plane $y = 0$, $\dot{y} < 0$, and prove the covering relations given below (the precise definition of the t-sets is given in Section 7):

$$N_0 \xRightarrow{H_1} M_0 \xRightarrow{H_2} N_0 \xRightarrow{H_1} M_1 \xRightarrow{H_2} N_1 \xRightarrow{H_1} M_2 \xRightarrow{H_2} N_2 \xRightarrow{H_1} M_3$$
$$\xRightarrow{H_2} N_3 \xRightarrow{H_1} M_4 \xRightarrow{H_2} N_4 \xRightarrow{H_1} M_5 \xRightarrow{H_2} N_5 \xRightarrow{H_1} M_6 \xRightarrow{H_2} N_5$$

$$L_3 \xRightarrow{H_2} K_0 \xRightarrow{H_1} L_3 \xRightarrow{H_2} N_3$$
$$K_1 \xRightarrow{H_1} L_1 \xRightarrow{H_2} K_2 \xRightarrow{H_1} L_2 \xRightarrow{H_2} N_3$$
$$K_3 \xRightarrow{H_1} L_4 \xRightarrow{H_2} K_3 \xRightarrow{H_1} L_5 \xRightarrow{H_2} N_3$$

But even to prove these covering relations with the half Poincaré maps turns out to be a difficult task, indeed some of the sets are quite close to the primaries (particularly the sets $L_4$ and $L_5$), and this causes the computation of rigorous bounds for the intersection of the trajectory with the Poincaré plane to be critical. To overcome this difficulty we have to use the alternative, but equivalent, definition of covering through the inverse of the Poincaré map or backcovering, see Definition 5.4. The covering relations we actually prove by computer assistance are as follows:

$$N_0 \xRightarrow{H_1} M_1 \xLeftarrow{H_2} N_1 \xRightarrow{H_1} M_2 \xLeftarrow{H_2} N_2 \xRightarrow{H_1} M_3 \xRightarrow{H_2} N_3 \xLeftarrow{H_1} M_4 \xRightarrow{H_2} N_4 \xLeftarrow{H_1} M_5 \xRightarrow{H_2} N_5$$

$$N_0 \xRightarrow{H_1} M_0, \qquad M_6 \xRightarrow{H_2} N_5, \qquad K_0 \xRightarrow{H_1} L_3 \xRightarrow{H_2} N_3 \tag{5.3}$$

$$K_1 \xRightarrow{H_1} L_1 \xLeftarrow{H_2} K_2 \xRightarrow{H_1} L_2 \xRightarrow{H_2} N_3, \qquad L_4 \xRightarrow{H_2} K_3 \xLeftarrow{H_1} L_5 \xRightarrow{H_2} N_3$$

Note that the coverings $L_3 \xLeftarrow{H_2} K_0$, $K_3 \xLeftarrow{H_1} L_4$, $M_0 \xLeftarrow{H_2} N_0$ and $N_5 \xLeftarrow{H_1} M_6$ follow from the coverings $K_0 \xRightarrow{H_1} L_3$, $L_4 \xRightarrow{H_2} K_3$, $N_0 \xRightarrow{H_1} M_0$ and $M_6 \xRightarrow{H_2} N_5$ respectively, by the symmetry of the t-sets and Lemma 5.7. As a result, we have to check 21 covering relations.

In order to have covering relations as described in Definition 5.3, we would need to check the image through the Poincaré map of the whole support of the t-sets. On the other hand, by Lemma 5.1 it is enough to perform rigorous computations on the boundaries of the t-sets, provided that we can prove that the map is defined on the whole supports. This is essential for the rigorous proof to be made in a reasonable time, since to check by computer assistance the definition of the Poincaré map on the supports of the t-sets would be very time consuming. In order to prove that the half Poincaré maps are defined on the supports of all t-sets we argue as follows. We only give the proof for the half Poincaré map $H_1$, the proof for $H_2$ being equivalent. First we have to prove that the trajectory does not collide with one primary. Then we observe that the projections of the trajectories we are interested in on the $(y, p_y)$ plane appear to be rotating around either a primary or the Lagrangian point $L_1$. If we can compute the angular velocity in the $(y, p_y)$ plane and prove that it is bounded away from 0 for long enough time, then it follows that the trajectory has to cross the Poincaré section eventually, hence the map $H_1$ is well defined. In [AZ] it was possible to prove by analytical computations that the projection on the $(x, p_x)$ plane of all trajectories have positive angular velocity with respect to the origin.

For the restricted 3-body problem such a general statement does not hold, therefore we have to proceed with a different method. First we compute the trajectory of the boundary of the t-set with rigorous bounds on the error. The trajectory of the whole boundary describes a "tube" in $\mathbb{R}^4$. More precisely, the intersection of the trajectory of the boundary with any hyperplane $\alpha p_y + \beta y = 0$ bounds a region on the energy surface; the union of all these intersections is what we call the tube. Since no trajectories can intersect, all trajectories starting in the interior of the support of the t-set have to remain inside the tube. If we can prove that the angular velocity in the interior of the tube is always bounded away from 0, say its absolute value is larger than $\varepsilon > 0$, it follows that the half Poincaré map is well defined on the whole interior of the support of the t-set; indeed any trajectory starting in the interior of the t-set must intersect the Poincaré plane in a time $T \leq \pi/\varepsilon$.

For a given point $(x, p_x) \in \partial N$ let

$$\Phi(x, p_x) = \varphi\left(x, p_x; [0, T_1]\right),$$

where $\varphi$ is the map defined in Section 3 and $T_1 = T_1(x, p_x)$ is the half Poincaré time also defined in Section 3. $\Phi(x, p_x)$ is the trajectory of the point $(x, p_x, 0, \sqrt{2\Omega(x,0) + h - p_x^2})$ under the flow induced by the equations from time 0 to the time $T_1$ when it reaches the Poincaré plane. If we can prove that $H_1(\partial N)$ exists, since the flows of two different points cannot intersect, unless they coincide, then $H_1(\partial N)$ bounds a region $M$ in the Poincaré plane. We want to prove that $H_1$ is defined on the whole set $|N|$ and $H_1(|N|) = M$. Let $\Sigma = |N| \cup \bigcup_{(x,p_x)\in\partial N} \Phi(x, p_x) \cup M$. The hypersurface $\Sigma \subset \mathbb{R}^4$ divides $\mathbb{R}^4$ into two connected regions; call $\Xi$ the bounded region. It is clear that for all $(x, p_x) \in |N| \setminus \partial N$ either there exists $H_1(x, p_x)$, or $\varphi\left(x, p_x, 0, \sqrt{2\Omega(x,0) + h - p_x^2}; T\right) \in \Xi$ for all $T \geq 0$. Let $\beta(x, p_x, y, p_y) = (2p_x + \Omega_y(x, y))y - p_y^2$. The next Lemma shows that, if $\Xi$ does not intersects the plane ($y = 0$, $p_y = 0$) and $|\beta(x, p_x, y, p_y)| \geq \delta > 0$ for all $(x, p_x, y, p_y) \in \Xi$, then the second case cannot happen.

**Lemma 5.8.** *Fix $h$ and $(x, p_x)$ such that $x \neq \{R_1, -R_2\}$ and $2\Omega(x,0) + h - p_x^2 > 0$; let $p_y = \sqrt{2\Omega(x,0) + h - p_x^2}$. Let $(x(T), p_x(t), y(T), p_y(T)) = \varphi(x, p_x, 0, p_y; T)$ be the flow induced by (1.1); let $\beta(T) = (2p_x(T) + \Omega_y(x(T), y(T)))y(T) - p_y^2(T)$. If there exists $\delta > 0$ such that $y^2(T) + p_y^2(T) > 0$ and $|\beta(T)| \geq \delta$ for all $T \in [0, \pi/\delta]$ (resp. $T \in [-\pi/\delta, 0]$), then the half Poincaré map $H_1$ (resp. $H_2$) in $(x, p_x)$ is well defined.*

*Proof.* Since $y^2(T) + p_y^2(T) > 0$ we can use polar coordinates. Let $\alpha(T)$ be the angle. Taking the derivative we have

$$\dot{\alpha} = \frac{\dot{p}_y y - p_y^2}{y^2 + p_y^2} = \frac{(2p_x + \Omega_y)y - p_y^2}{y^2 + p_y^2}.$$

Since $y^2 + p_y^2$ is bounded and $(2p_x + \Omega_y)y - p_y^2$ is bounded away from 0, then for some $T$ we have $|\alpha(T) - \alpha(0)| \geq \pi$ and the map is defined. □

**Lemma 5.9.** *For all covering relations listed in (5.3) the half Poincaré maps are defined on all the supports of all t-sets.*

*Proof.* Consider a covering relation with the $H_1$ map, the $H_2$ map and the back-covering relations being equivalent. The proof is by computer assistance and it is performed with the following procedure. Choose a t-set and let $(x, p_x)$ the middle point of its support. Let $\bar{\varphi}(x, p_x, t)$ be the approximate flow, computed by some suitable algorithm (we used Mathematica 4.0 for this purpose). Let $\bar{T}_1$ be the approximate Poincaré time, i.e. $\bar{\varphi}_3(x, p_x, t) > 0$ for all $0 < t < \bar{T}_1$ and $\bar{\varphi}_3(x, p_x, \bar{T}_1) = 0$. For an integer $N > 0$ let $t_i = i\bar{T}_1/N$, $i = 0, N$,

$$\bar{\Xi}_i = \left\{(x, p_x, y, p_y) \in \mathbb{R}^4 : d((x, p_x, y, p_y), \bar{\varphi}(x, p_x, t_i)) \leq 1\right\},$$

where $d$ is a distance in $\mathbb{R}^4$ defined by

$$d((x^1, p_x^1, y^1, p_y^1), (x^2, p_x^2, y^2, p_y^2)) = \max(a_1|x^1 - x^2|, a_2|p_x^1 - p_x^2|, a_3|y^1 - y^2|, a_4|p_y^1 - p_y^2|)$$

and $\{a_i\}$ are positive constants. Let

$$\bar{\Xi} = \bigcup_{i=0,N} \bar{\Xi}_i.$$

$\bar{\Xi}$ is a very rough approximation of the flow of the support of the t-set. If all the constants entering in the definition of $\bar{\Xi}$ are chosen appropriately, one can expect that the true trajectory of the support of the t-set is entirely contained in $\bar{\Xi}$. With computer assistance we checked rigorously that the true trajectory of the boundary of the t-sets never leaves $\bar{\Xi}$. Using interval arithmetics algorithms we compute

$$\delta_i = \max_{(x,p_x,y,p_y)\in\bar{\Xi}_i} (2p_x + \Omega_y(x, y))y - p_y^2.$$

By computer assistance we prove that $\max_i \delta_i$ is strictly negative for all t-sets. The proof is complete by Lemma 5.8. □

The previous lemmas yield the following:

**Theorem 5.10.** *All covering relations listed in (5.3) hold.*

*Proof.* We check with computer assistance that conditions **a′** and **b** in Lemma 5.1 or conditions **c** and **d** in Definition 5.4 are verified for all the covering relations listed in (5.3). The proof is completed by Lemma 5.9. □

Consider the collection of the t-sets $N_0$, $N_1$, $N_2$, $N_3$, $N_4$, $N_5$, $\tilde{N}_1$, $\tilde{N}_2$, $\tilde{N}_3$, $\tilde{N}_4$, $K_0$, $K_1$, $K_2$, $K_3$ in this order and let

$$T = \begin{bmatrix}
1 & 1 & 0 & 0 & 0 & 0 & 0 & 0 & 0 & 0 & 0 & 0 & 0 & 0 \\
0 & 0 & 1 & 0 & 0 & 0 & 0 & 0 & 0 & 0 & 0 & 0 & 0 & 0 \\
0 & 0 & 0 & 1 & 0 & 0 & 0 & 0 & 0 & 0 & 0 & 0 & 0 & 0 \\
0 & 0 & 0 & 0 & 1 & 0 & 0 & 0 & 0 & 0 & 0 & 0 & 0 & 0 \\
0 & 0 & 0 & 0 & 0 & 1 & 0 & 0 & 0 & 0 & 0 & 0 & 0 & 0 \\
0 & 0 & 0 & 0 & 0 & 1 & 0 & 0 & 0 & 1 & 0 & 0 & 0 & 0 \\
1 & 0 & 0 & 0 & 0 & 0 & 0 & 0 & 0 & 0 & 0 & 0 & 0 & 0 \\
0 & 0 & 0 & 0 & 0 & 0 & 1 & 0 & 0 & 0 & 0 & 0 & 0 & 0 \\
0 & 0 & 0 & 0 & 0 & 0 & 0 & 1 & 0 & 0 & 1 & 0 & 1 & 1 \\
0 & 0 & 0 & 0 & 0 & 0 & 0 & 0 & 1 & 0 & 0 & 0 & 0 & 0 \\
0 & 0 & 0 & 1 & 0 & 0 & 0 & 0 & 0 & 1 & 0 & 0 & 0 & 0 \\
0 & 0 & 0 & 0 & 0 & 0 & 0 & 0 & 0 & 0 & 0 & 0 & 1 & 0 \\
0 & 0 & 0 & 1 & 0 & 0 & 0 & 0 & 0 & 0 & 0 & 1 & 0 & 0 \\
0 & 0 & 0 & 1 & 0 & 0 & 0 & 0 & 0 & 0 & 0 & 0 & 0 & 1
\end{bmatrix}.$$

Theorem 5.10 and the definition of generic covering yield the following

**Corollary 5.11.** *All relations in (5.1) and (5.1) hold with generic covering instead of covering or backcovering and $T$ is the transition matrix associated with the map $P$.*

5.4. **Symbolic dynamics.** Let $\Sigma_A$ be the set of all admissible sequences of fourteen symbols ($N_0$, $N_1$, $N_2$, $N_3$, $N_4$, $N_5$, $\tilde{N}_1$, $\tilde{N}_2$, $\tilde{N}_3$, $\tilde{N}_4$, $K_0$, $K_1$, $K_2$, $K_3$) with respect to $T$ (recall Definition 5.8). Corollary 5.11 yields the following result on symbolic dynamics:

**Theorem 5.12.** *For any biinfinite sequence $\{x_n\} \in \Sigma_A$ there exists an orbit of (1.1) which crosses the Poincaré plane in the t-sets $N_0, \ldots, N_5$, $\tilde{N}_1, \ldots, \tilde{N}_4$ and $K_0, \ldots, K_3$ in the order prescribed by the sequence. Furthermore, if $\{x_n\}$ is periodic of period $m$, such orbit is periodic as well and its trajectory on the Poincaré section also has period $m$.*

*Proof.* The supports of the t-sets $N_0, \ldots, N_5$, $\tilde{N}_1, \ldots, \tilde{N}_4$ and $K_0, \ldots, K_3$ are disjoint; let $N$ be the union of all such supports and let $\pi : \mathrm{Inv}(N, P) \to \Sigma_A$ be the projection defined as in Definition 5.9. The result follows by Theorems 5.6 and Corollary 5.11. □

**Corollary 5.13.** *Equation (1.1) admits infinitely many periodic solutions.*

This result is important since it gives the complete picture of the chaotic behavior of the system and it yields the result on topological entropy we present in the next subsection. On the other hand, a better visual image of the result can be obtained by the following observations. Recall that the sets $N_0$, $N_5$, $K_0$ and $K_3$ $P-$cover themselves and $K_1$ $P^2$-covers itself. Hence there exists (at least) a periodic orbit for each of these sets which crosses the hyperplane $y = 0$ with positive $y$ velocity in the set (or alternatively in the sets $K_1$ and $K_2$). The orbit crossing $N_0$ is the orbit symmetric to $U_2$ with respect to the $y$ axis, while the orbits crossing $N_5$, $K_0$, $K_1$ and $K_3$ are $U_2$, $L$, $D_1$ and $U_1$ respectively. We point out that without considering the derivative of the Poincaré map we cannot exclude that other periodic orbits cross those t-sets. It would be indeed rather easy to rigorously compute the derivative of the Poincaré map, but this is not relevant for the discussion that follows and such computation will be treated in another paper.

From the transition matrix it also follows that $N_0 \overset{P^5}{\Longleftrightarrow} N_5 \overset{P^5}{\Longleftrightarrow} N_0$ and $K_i \overset{P^3}{\Longleftrightarrow} N_5 \overset{P^3}{\Longleftrightarrow} K_i$, $i = 0, 2, 3$, therefore if we consider the t-sets $N_0$, $N_5$, $K_0$, $K_2$ and $K_3$ the transition matrix with respect to the map $P^{30}$ is a $5 \times 5$ matrix with all entries equal to 1, that is each of these t-sets $P^{30}-$covers itself and all the other t-sets. Using again Theorem 5.6 we infer that for all biinfinite sequences of five symbols there exists a solution of equation (1.1) which crosses the t-sets in the same order. In other words, the system admits solutions which come close to any of these orbits in any prescribed order. We point out that this result concern the thirtieth interate of the Poincaré map, therefore only one crossing of the Poincaré section every 30 should be considered.

Finally, we point out that all these results concern the existence of orbits, not uniqueness.

5.5. **Topological entropy estimates.** The results presented in the previous subsection yield a lower estimate of the topological entropy of the map $P$ by the following lemma:

**Lemma 5.14.** *Let $f : X \to X$ be a continuous map. Let $S \subset X$ be an invariant set, let $A$ be an $n \times n$ matrix such that there exists a surjective map $\pi : S \to \Sigma_A$ satisfying $\sigma \circ \pi = \pi \circ f$. Then the topological entropy of $f$ is larger than $\ln(\max\{|\lambda_i|, \quad \lambda_i$ is an eigenvalue of $A\})$.*

*Proof.* The proof is an easy consequence of Theorem 7.13 in [WA]. □

**Corollary 5.15.** *The topological entropy of $P$ is larger than* 1.62746.

*Proof.* By Lemma 5.14 a lower bound for $h_t$ is given by the maximal norm eigenvalue of the transition matrix $T$. The characteristic polynomial is $p(x) = 1 - 2\,x + 2\,x^3 + x^4 - 3\,x^5 - 3\,x^6 + 7\,x^7 - 4\,x^8 + 4\,x^9 - 5\,x^{10} + 5\,x^{12} - 4\,x^{13} + x^{14}$, and since $p(1.62747) > 0$ and $p(1.62746) < 0$, then $h_t(P) > 1.62746$. □

## 6. Further developments

Although the 3-body problem is very old, there are still many open problems. Computer assisted techniques can give rigorous proofs to some conjectures. This paper is focused on results which require only $C^0$ computations. We think that a much richer symbolic dynamics can be proved with these methods at different energy levels: this paper should only be considered as an example, although in the author opinion Theorem 5.12 and Corollary 5.15 have a value in themselves. The rigorous computation of the derivatives of the Poincaré map can provide information on the linear stability of periodic orbits and yield existence results for homoclinic and heteroclinic orbits. We will treat these topics in a forthcoming paper, together with the dependence of the results presented in this paper to the ratio of the masses of the primaries. A further development under study concerns the study of bifurcations with the energy as a parameter value.

## 7. Computational details

7.1. **Description of the t-sets.** The procedure used to build all t-sets is described in [AZ] for the Hénon-Heiles Hamiltonian and we refer to that paper for the details. Here we only mention that the t-sets $N_i$ and $K_i$ are either centered on the fixed or period 2 points of the Poincaré map, or placed along the invariant manifolds of such points with the sides approximately parallel to the manifolds. The sets $M_i$ and $L_i$ instead are built "by hand" trying to interpolate the image of the other sets. We recall that the Poincaré section is not connected, since the lines $x = \pm 2^{-2/3}$ corresponding to the position of the primaries must be excluded. It follows that the invariant manifolds of points belonging to different connected components cannot cross. This fact does not influence the topological method employed here, since a t-set and its image may be on different connected components, bypassing the fact that the Poincaré section passes through the primaries. The actual t-sets used in the proofs are as follows:

**Definition 7.1.** *The triple sets are defined by giving the coordinates of the center* $(x, y)$*, the length of the sides* $(l_x, l_y)$ *and the angular coefficients of the sides* $(\alpha, \beta)$ *as follows:*

$N_0$: $(x, y) = (-0.7942, 0.)$, $(l_x, l_y) = (.007, .007)$, $(\alpha, \beta) = (.9325, -.9325)$.
$N_1$: $(x, y) = (-.7849, .0129)$, $(l_x, l_y) = (.002, .005)$, $(\alpha, \beta) = (-2.136, 2.136)$.
$N_2$: $(x, y) = (-.7589, .05086)$, $(l_x, l_y) = (.002, .005)$, $(\alpha, \beta) = (-1.812, 1.981)$.
$N_3$: $(x, y) = (.08, .125)$, $(l_x, l_y) = (.01, .04)$, $(\alpha, \beta) = (1.124, 2.009)$.
$N_4$: $(x, y) = (.131, .028)$, $(l_x, l_y) = (.0076, .012)$, $(\alpha, \beta) = (-2.09, 2.11)$.
$N_5$: $(x, y) = (.1471, 0.)$, $(l_x, l_y) = (.015, .015)$, $(\alpha, \beta) = (-2.07, 2.07)$.
$K_0$: $(x, y) = (.04873, 0.)$, $(l_x, l_y) = (.0002, .0002)$, $(\alpha, \beta) = (1.818, -1.818)$.
$K_1$: $(x, y) = (.04775, 0.)$, $(l_x, l_y) = (.0009, .0009)$, $(\alpha, \beta) = (1.131, -1.131)$.
$K_2$: $(x, y) = (-.7365, 0.)$, $(l_x, l_y) = (.0012, .0012)$, $(\alpha, \beta) = (-2, 2)$.
$K_3$: $(x, y) = (.06712, 0.)$, $(l_x, l_y) = (.01, .01)$, $(\alpha, \beta) = (1.835, -1.835)$.
$M_0$: $(x, y) = (-.1471, 0.)$, $(l_x, l_y) = (.0065, .0065)$, $(\alpha, \beta) = (1, -1)$.
$M_1$: $(x, y) = (-.1424, .007289)$, $(l_x, l_y) = (.0026, .0065)$, $(\alpha, \beta) = (.9, -1.05)$.
$M_2$: $(x, y) = (-.1311, .0279)$, $(l_x, l_y) = (.0025, .005)$, $(\alpha, \beta) = (4.2, 2)$.

*$M_3$: $(x,y)=(-.0792,.1267)$, $(l_x,l_y)=(.002,.005)$, $(\alpha,\beta)=(4.2,1.8)$.*
*$M_4$: $(x,y)=(.7592,.05021)$, $(l_x,l_y)=(.006,.018)$, $(\alpha,\beta)=(1.2,2.2)$.*
*$M_5$: $(x,y)=(.7854,.01406)$, $(l_x,l_y)=(.01,.008)$, $(\alpha,\beta)=(4.2,2.2)$.*
*$M_6$: $(x,y)=(0.7942,0.)$, $(l_x,l_y)=(.013,.013)$, $(\alpha,\beta)=(1,-1)$.*
*$L_1$: $(x,y)=(-.0742,-.1015)$, $(l_x,l_y)=(.0009,.0009)$, $(\alpha,\beta)=(1.131,-1.131)$.*
*$L_2$: $(x,y)=(-0.0738,.1025)$, $(l_x,l_y)=(.002,.002)$, $(\alpha,\beta)=(-1.85,1.85)$.*
*$L_3$: $(x,y)=(-.04873,0.)$, $(l_x,l_y)=(.01,.01)$, $(\alpha,\beta)=(1.819,-1.819)$.*
*$L_4$: $(x,y)=(.6563,0.)$, $(l_x,l_y)=(.005,.005)$, $(\alpha,\beta)=(-1.49,1.49)$.*
*$L_5$: $(x,y)=(.6545,.022)$, $(l_x,l_y)=(.01,.01)$, $(\alpha,\beta)=(-1.49,1.49)$.*

*The left (resp. right) edge of each set is the segments whose end points are $(x+l_x\cos\alpha+l_y\cos\beta, y+l_x\sin\alpha+l_y\sin\beta)$ and $(x+l_x\cos\alpha-l_y\cos\beta, y+l_x\sin\alpha-l_y\sin\beta)$ (resp. $(x-l_x\cos\alpha+l_y\cos\beta, y-l_x\sin\alpha+l_y\sin\beta)$ and $(x-l_x\cos\alpha-l_y\cos\beta, y-l_x\sin\alpha-l_y\sin\beta)$). The boundaries of the left and right sides of each t-set are the lines crossing two opposite vertices of the support, see Figure 3.*

*Furthermore, we denote $\tilde N_i$, $i=1,\dots 4$, the sets obtained by reflecting the sets $N_i$ with respect to the $x$-axis, see Definition 5.2. Note that the sets $N_0$, $N_5$, $K_0$, $K_1$, $K_2$, $K_3$, $K_3$ are by definition symmetric with respect to the $x$ axis.*

7.2. **Computation techniques.** We describe here the algorithm used in the computer assisted proof of Theorem 5.10. The proof of Theorem 4.1 is equivalent. By Lemma 5.1 the proof consists in checking that the images through the Poincaré map of the edges of the t-sets lie in some assigned regions of the plane as described by the covering relations. We do not know the exact images of such edges, since no analytical solution of the equation is available, therefore all we can do is to estimate the trajectories and compute the intersections with the Poincaré plane with rigorous error bounds. In order to compute the image of a side, we partition it in segments with small enough length and we check that every such segment is mapped in the correct region. To compute the image of a small segment, we enclose it in an interval set (a rectangle) and we compute its trajectory with a Taylor-Lohner algorithm using interval arithmetics. More precisely, we start with a Taylor method of order 12, i.e. we estimate the trajectory of an interval by using the Taylor expansion of order 12 and we estimate the error by the Lagrange remainder. If $h$ is the time step, we compute a rough but rigorous enclosure $D$ of the trajectory at times $[0,h]$, that is an interval set $D$ such that the solution of the equation lies in $D$ for all times between 0 and $h$, and by Lagrange theorem we estimate the error we make neglecting the remaining terms of the Taylor expansion by computing $x^{(13)}(D)\frac{h^{13}}{13!}$, where $x^{(13)}(D)$ (which is an interval enclosing all possible values assumed by the 13th derivative of the trajectory, therefore enclosing the Lagrange remainder) is computed using a recursive algorithm for the time derivatives of the solutions.

The interval arithmetics algorithms address the problem of computing the trajectory of a segment and of keeping track of the errors in an elegant and rigorous way, but they introduce another problem. Indeed, even in the simplest dynamical system, the procedure described above leads to a very rough estimate of trajectories, due to the wrapping effect which makes the bounds on the error grow exponentially. The problem has been strongly contained by introducing the half Poincaré maps and the backcovering relation, see the discussion in Section 5, but these techniques do not suffice. We obtain another significant reduction of the wrapping effect by using the Lohner algorithm. We refer to see Section 6 in [AZ] and references cited there for a discussion of interval arithmetics and wrapping effect and for a description of the Lohner algorithm employed. See also [MZ] for a discussion on interval arithmetics and [L] for the Lohner algorithm. We point out that using the Lohner algorithm to compute directly the Poincaré map, instead of the half Poincaré maps or their inverse, the computation time to perform all the proofs is unrealistic on current desk computers, therefore the definitions of the half Poincaré maps and of backcovering are essential to perform the proofs.

The round-off errors are taken care directly by suitable C++ libraries and by Mathematica. Such errors may vary by changing computer and/or operationg system, but since they are usually very small when compared to the wrapping effect and since all proofs go through with a relatively large margin, we expect that the proof can be easily reproduced on any recent computer.

The typical time step used in the computation of the images of the t-sets is $dt = 10^{-2}$, but in few cases we had to use some lower value, down to $dt = 5.10^{-4}$. Each side of the t-sets has been usually divided in 100 to 2000 segments, depending on the apparent value of the Lipschitz constant of the map in the area considered. In a few difficult cases we had to divide each side of a t-set in 5000 segments. Most of the computations used to obtain the bounds for the location of the orbits in Table 1 used $dt = 2.10^{-5}$.

To perform the proofs the author implemented a version of the whole algorithms in a combination of Mathematica and C++ under the Linux O.S. More precisely, Mathematica has been only used to handle all the data and to perform a few algorithms which are less demanding for the CPU, but more complicated to implement. Furthermore Mathematica has been used to make all numerical experiments used to choose the t-sets and to draw the pictures. On the other hand C++ has been used for the heavy interval arithmetic computations, where it offered a much better speed. The connection between the two languages is obtained by MathLink. We wish to point out that the full proof took almost a month of CPU time on a machine equipped with a 1GHz Pentium III processor. In fact a large amount of the time is used for the computations involving the set $L_5$ whose trajectory is very close to one primary. A full C++ algorithm would reduce the time at the price of a much more complicated and less user-friendly programming. We think that this sharing of tasks is almost optimal, as far as computational speed and simplicity of programming and data handling is concerned. The reader who desires to reproduce the computer assisted proofs in this paper without writing the program can use the Mathematica notebook which is attached to this preprint. The notebook is provided with comments and instructions. By using the commands in the notebook it is possible to make an independent computation of the images of the sides of the t-sets and an independent check of Theorem 4.1.

**Acknowledgement.** The author is very grateful to P. Zgliczyński for many discussions and for the interval arithmetic C++ libraries.

Dipartimento di Scienze e T.A., C.so Borsalino 54, 15100 Alessandria Italy
*E-mail address*: `gianni@unipmn.it` *Web page:* `http://www.mfn.unipmn.it/~gianni/eng.html`